\documentclass[12pt]{article}
\usepackage{hyperref}
\usepackage{amsmath}
\usepackage{ntheorem}
\usepackage{mathtools}
\usepackage{epsf}
\usepackage{theorem}
\usepackage{indentfirst}
\usepackage{latexsym}
\usepackage{amssymb}
\usepackage[dvips]{graphicx}
\usepackage{flafter}
\usepackage{caption}
\usepackage{textcomp}
\usepackage{supertabular}
\usepackage{longtable, booktabs}
\usepackage{hhline}
\usepackage{bigstrut,bigdelim}
\usepackage{rotating}
\usepackage{multicol}
\usepackage{multirow}
\usepackage{makeidx}
\usepackage{epsfig}
\usepackage{fancyhdr}
\usepackage{threeparttable}
\usepackage{longtable}
\usepackage{mathrsfs}
\usepackage[all]{xy}
\usepackage{xcolor}
\newif\ifger
\gerfalse
\newtheorem{definition}{Definition}[section]
\newtheorem{theorem}{Theorem}[section]

\newtheorem{corollary}[theorem]{Corollary}
\newtheorem{remark}{Remark}[section]

\newtheorem{proposition}[definition]{Proposition}

\newtheorem*{availability of data and material}{Availability of Data and Material}
\newtheorem*{conflict of interest}{Conflict of interest}
\begin{document}
\baselineskip=19pt

\title{Some permutation polynomials via linear translators}

\author{Xuan Pang$^a$, Pingzhi Yuan$^a${\footnote{Corresponding author. E-mail: yuanpz@scnu.edu.cn(P. Yuan)}, Hongjian Li$^b$.} \\
\small \it  $^a$School of Mathematical Sciences, South China Normal University, \\
\small \it  Guangzhou, 510631, P. R. China\\
\small \it  $^b$School of Mathematics and Statistics, Guangdong University of Foreign Studies, \\
\small \it Guangzhou, 510006, P. R. China\\}

\date{}
\maketitle

\begin{abstract}

Permutation polynomials with explicit constructions over finite fields have long been a topic of great interest
in number theory. In recent years, by applying linear translators of functions from $\mathbb{F}_{q^n}$ to $\mathbb{F}_q$,
many scholars constructed some classes of permutation polynomials. Motivated by previous works, we first naturally extend
 the notion of linear translators and then construct some permutation polynomials.



\end{abstract}

\section{Introduction}

   Let $\mathbb{F}_q$ be the finite field with $q$ elements and $\mathbb{F}_q^*$ be its multiplicative group consisting of all nonzero elements of $\mathbb{F}_q$, where $q$ is a prime power. A polynomial $f(x)\in \mathbb{F}_q[x]$ is called a {\it permutation polynomial} (PP) of $\mathbb{F}_q$ if its associated polynomial mapping 
    \begin{align*}
    f:~ &\mathbb{F}_q \longrightarrow \mathbb{F}_q\\
    & ~\alpha \longmapsto f(\alpha)
   \end{align*}
   induces a permutation. It is widely known that for any permutation polynomial over $\mathbb{F}_q$, there exists a unique polynomial denoted by $f^{-1}(x)$ over $\mathbb{F}_q$ such that $f\circ f^{-1}(x)=f^{-1}\circ f(x)=x$, in the sense of reducing modulo $x^q-x$, and it is said to be the compositional inverse of $f(x)$. In particular, we call $f(x)$  an involution if its compositional inverse is itself, i.e., $f(x)=f^{-1}(x)$.

   Constructing permutation polynomials over finite fields is an ancient and significant topic due to their broad applications, especially in the fields of  coding theory
    \cite{code1-Harish,code2-Laigle-Chapuy,code3-Li}, cryptography \cite{cryp1-Niederreiter,cryp2-Schwenk,cryp3-Singh} and
    so on. However, it is still a challenge to determinate or construct some new permutation polynomials,
    despite the efforts of numerous scholars. Indeed, the study of permutation polynomials is usually considered to have
    begun in the 19th century when Dickson \cite{Lidl} applied Hermite's criteria to derive all normalized permutation
     polynomials of degree at most 5. In 2008, Zieve \cite{Zieve-criteria} provided yet another effective technique
     for constructing new classes of permutation polynomials. Motivated later by Marcos \cite{Marco}
     and Zieve's \cite{Zieve2010} work, Akbary, Ghioca and Wang \cite{AGW} proposed a more general method named AGW criterion,
      which has led to many new results in recent years. More recently, Yuan \cite{yuan-local} also presented a local method
      to obtain permutation polynomials as well as their compositional inverses.

   Among the approaches for constructing permutation polynomials, explicit constructions are of a great interest.
   This brings us to the concept of linear translators which originated from the Boolean function, and then was introduced
    over finite fields by Charpin and Kyureghyan. We say that $\gamma\in \mathbb{F}_{q^n}^*$ is a $b$-linear translator of
    the function $f$: $ \mathbb{F}_{q^n} \rightarrow \mathbb{F}_q$ if $f(x+u\gamma)-f(x)=ub$ holds,
    for all $x\in \mathbb{F}_{q^n}$, $u\in \mathbb{F}_q$ and $b$ fixed in $\mathbb{F}_q$.
    In 2008, Charpin and Kyureghyan \cite{Ch-Ky2008} studied permutation polynomials of the shape $G(x)+\gamma Tr(H(x))$ over finite fields with even characteristic by mean of the linear structures, 
   and obtained six classes of such permutations. And the two authors \cite{Ch-Ky2009} in turn generalized this work
   to the situation in which  finite fields with odd characteristic. Meanwhile, in \cite{ky2011},
   Kyureghyan also construcetd large families of permutation polynomials of the form $L(x)+L(\gamma)h(f(x))$,
   where $L(x)=\sum\limits_{i=0}\limits^{n-1}{a_ix^{q^i}}$ is an $\mathbb{F}_q$-linear permutation over $\mathbb{F}_{q^n}$.

 \begin{theorem}\cite[Theorem 1]{ky2011}\label{cy-th}
 Let $f:\mathbb{F}_{q^n}\rightarrow\mathbb{F}_q$, $h:\mathbb{F}_q\rightarrow\mathbb{F}_q$ be arbitrary mapping. Let $L:\mathbb{F}_{q^n}\rightarrow\mathbb{F}_{q^n}$ be an $\mathbb{F}_{q}$-linear permutation of $\mathbb{F}_{q^n}$. If $b\in \mathbb{F}_q$, and $\gamma\in \mathbb{F}_{q^n}^*$ is a $b$-linear translator of f, then
 $$G(x)=L(x)+L(\gamma)h(f(x))$$
 permutes $\mathbb{F}_{q^n}$ if and only of $g(x)=x+bh(x)$ permutes $\mathbb{F}_{q}$.
 \end{theorem}
 Since then, a large number of new bent functions and permutation polynomials have been constructed using a similar way
   by many researchers \cite{Cepak2017,Mesnager2017,qin2015,qin2020,xie2021}.

  Recently, considering the case when $f(x+\gamma u)-f(x)=u^{p^i}b$,  Cepak, Pasalic and
  Moratovi$\acute{\rm {c}}$-Ribi$\acute{\rm{c}}$ \cite{Cepak2019} generalized the definition of linear translators
  to Frobenius translators, and gave some constructions of permutation polynomials. Moreover, Mesnager, Kaytanc{\i}, $\ddot{\rm O}$zbdak
  further extended the notion of linear translators in \cite{Mesnager2019}. Motivated by them, we make a more natural generalization of the concept of linear translators, and then obtain a few permutation polynomials
  over finite fields.

The structure of this paper is organized as follows. In Section 2, we introduces the basic definition of linear translators in our paper, as well as some of their properties. In Section 3, we present rather more general versions of permutation polynomials from \cite{Cepak2019} and \cite{qin2020}.

\section{Preliminaries}

 We begin this section by providing the following definition of linear translators, which is somewhat different from  \cite{Ch-Ky2009,Cepak2019,Mesnager2019}.

 \begin{definition}\label{def}
 Let $f: \mathbb{F}_{q^n} \rightarrow \mathbb{F}_q$, and $A(x)$ be an additive permutation polynomial over $\mathbb{F}_q$. A non-zero element $\gamma\in \mathbb{F}_{q^n}$ is called a $(b,A)$-linear translator of the function $f$ if $$f(x+u\gamma)-f(x)=bA(u)$$
 holds for any $x\in\mathbb{F}_{q^n}, u\in \mathbb{F}_q$ and a fixed $b\in \mathbb{F}_q$.
 \end{definition}

\begin{remark}\label{Frobenius translator}
The notations and assumptions are the same as defined above. Notice that taking $A(x)=x^{p^i}$  agrees with Definition 2 in \cite{Cepak2019}, that is , $\gamma$ is an $(i,b)$-Frobenius translator of $f$, where $p$ is the characteristic of $\mathbb{F}_q$. In addition, we assert that $f$ must be surjective if there exists a $(b,A)$-linear translator of $f$, since $A(x)$ is a permutation polynomial over $\mathbb{F}_q$.
\end{remark}

In the following, we state some basic results  which will be used in this paper.

\begin{proposition}\label{pro2.2}
Let $f$ be a function from $\mathbb{F}_{q^n}$ to $\mathbb{F}_q$, and let $\gamma_1, \gamma_2 \in \mathbb{F}_{q^n}^*$ with $\gamma_1+\gamma_2\neq 0$. If $\gamma_1$ is a $(b_1,A)$-linear translator of $f$, and $\gamma_2$ is a $(b_2,A)$-linear translator of $f$, where $A(x)\in \mathbb{F}_q[x]$ is an additive permutation polynomial, then
$\gamma_1+\gamma_2$ is a $(b_1+b_2,A)$-linear translator of $f$, and for any $c\in F_p^*$, $c\gamma_1$ is a $(cb_1,A)$-linear translator of $f$.
\end{proposition}
\textbf{Proof.}\quad Let $\gamma_i$ be a $(b_i,A)$-linear translator of $f:\mathbb{F}_{q^n}\rightarrow \mathbb{F}_q$, where $i\in \{1,2\}$. Then for all $x\in \mathbb{F}_{q^n}, u\in\mathbb{F}_q$, we have
\begin{align*}
f(x+u(\gamma_1+\gamma_2))&=f(x+u\gamma_1+u\gamma_2)-f(x+u\gamma_1)+f(x+u\gamma_1)\\
   &=b_2A(u)+f(x+u\gamma_1)-f(x)+f(x)\\
   &=f(x)+b_1A(u)+b_2A(u)\\
   &=f(x)+(b_1+b_2)A(u).
\end{align*}
 This implies that $\gamma_1+\gamma_2$ is a $(b_1+b_2,A)$-linear translator of $f$. It remains to show that for an arbitrary non-zero element $c\in \mathbb{F}_p$, $c\gamma_1$ is also a linear translator of $f$. Giving any $c$ in $\mathbb{F}_p^*$, we see that
 \begin{align*}
 f(x+u(c\gamma_1))=f(x)+b_1A(uc)=f(x)+b_1cA(u),
  \end{align*}
since $\gamma_1$ is a linear translator of $f$. Thus, we have done.   $\hfill\square$

\begin{remark}
The set of all $(b,A)$-linear translators of $f$, along with zero element, forms an $\mathbb{F}_{p}$-linear subspace.
\end{remark}

Based on the notion in Definition \ref{def}, we have a generally consequence of Theorem \ref{cy-th} as below. 
In fact, this result can be obtained in a similar argument to \cite[Proposition 1]{Cepak2019} or Theorem \ref{cy-th}, so we omit its proof in the following.
\begin{theorem}
Let $L$ be an $\mathbb{F}_q$-linear permutation of $\mathbb{F}_{q^n}$, $A(x)$ an additive permutation of $\mathbb{F}_q$. Let $f:\mathbb{F}_{q^n}\rightarrow\mathbb{F}_{q}$, $h:\mathbb{F}_q\rightarrow\mathbb{F}_q$, $\gamma\in\mathbb{F}_{q^n}^*$ and $b$ fixed in $\mathbb{F}_{q}$. If $\gamma$ is a $(b,A)$-linear translator of $f$, then the mapping
$$G(x)=L(x)+L(\gamma)A^{-1}(h(f(x)))$$
permutes $\mathbb{F}_{q^n}$ if and only if $g(x)=x+bh(x)$ permutes $\mathbb{F}_q$, where $A^{-1}(x)$ is the compositional inverse of $A(x)$.
\end{theorem}

\begin{proposition}\label{pro2.3}
Let $f$ be a function from $\mathbb{F}_{q^n}$ to $\mathbb{F}_q$ and $A(x)$  an additive permutation polynomial
over $\mathbb{F}_{q}$. Let $\gamma_1,\gamma_2,\cdots \gamma_m\in \mathbb{F}_{q^n}^*$,
and $b_1,b_2,\cdots b_m\in \mathbb{F}_q$. If $\gamma_i$ is a $(b_i,A)$-linear translator of $f$,
where $i\in \{1,2,\cdots, m\}$, then for any $x\in\mathbb{F}_{q^n}, u_i\in \mathbb{F}_q$,
$$f(x+\sum\limits_{i=1}\limits^{m}{u_i\gamma_i})-f(x)=\sum\limits_{i=1}\limits^{m}{b_iA(u_i)}.$$
\end{proposition}
\textbf{Proof.}\quad Repeatedly applying Proposition \ref{pro2.2} yields this conclusion. $\hfill\square$

\medskip


Notice that the trace function $Tr^n_m(x)$ from $\mathbb{F}_{p^n}$ to $\mathbb{F}_{p^m}$, where $m|n$, is defined by $$Tr^n_m(x)=x+x^{p^m}+x^{p^{2m}}+\cdots+x^{p^{(\frac{n}{m}-1)m}}, ~~x\in \mathbb{F}_{p^n}.$$
For $m=1$, we abbreviate and write $Tr(x)=Tr^n_1(x)$, called the absolute trace.

\begin{proposition}\label{example}
For positive integers $m,n,s$, with $m|n, 1\leq s\leq m-1$, let $f:\mathbb{F}_{p^n}\rightarrow \mathbb{F}_{p^m}$ be a mapping defined by $f(x)=Tr_m^n(x^{p^s}-\alpha\gamma^{p^s-1}x)$, where $\gamma\in \mathbb{F}_{p^n}^*$ and $\alpha\in \mathbb{F}_{p^m}^*$ is
 not a $(p^s-1)$-th power, that is, $\alpha \neq \beta^{p^s-1}$ for arbitrary $\beta\in\mathbb{F}_{p^m}$. Then $\gamma$ is a $(Tr_m^n(\gamma^{p^s}), x^{p^s}-\alpha x)$-linear translator of $f$.
\end{proposition}
\textbf{Proof.}\quad Let $\gamma\in\mathbb{F}_{p^n}^*, u\in\mathbb{F}_{p^m}$. Then
\begin{align*}
f(x+u\gamma)-f(x)&=Tr_m^n((x+u\gamma)^{p^s}-\alpha\gamma^{p^s-1}(x+u\gamma))-Tr_m^n(x^{p^s}-\alpha\gamma^{p^s-1}x)\\
 &=Tr_m^n(u^{p^s}\gamma^{p^s}-\alpha u\gamma^{p^s})\\
 &=(u^{p^s}-\alpha u)Tr_m^n(\gamma^{p^s}).
\end{align*}
Define $A(u)=u^{p^s}-\alpha u$. We easily see that $A(u)$ is an additive permutation of $\mathbb{F}_{p^m}$, since $\xi_1^{p^s}-\alpha \xi_1 \neq\xi_2^{p^s}-\alpha\xi_2$ for any $\xi_1\neq \xi_2$, otherwise $\alpha=(\xi_1-\xi_2)^{p^s-1}$,  which is a contradiction. Therefore, for any $x\in\mathbb{F}_{p^n}$ and any $u\in\mathbb{F}_{p^m}$, we have $f(x+u\gamma)-f(x)=Tr_m^n(\gamma^{p^s})A(u)$. This concludes the proof.
 $\hfill\square$


\section{Permutation polynomials via linear translators}

  As a matter of notation $A(x)$ or $A_i(x)$ will be an additive permutation polynomial over $\mathbb{F}_q$ throughout this paper.
  In this section, we first gives a generalized result for permutations of $\mathbb{F}_{2^n}$ in \cite{Cepak2019} and \cite{Mesnager2017}, as well as their induced bent functions, using $(b,A)$-linear translators of $f:\mathbb{F}_{2^n}\rightarrow \mathbb{F}_{2^m}$, where $q=2^m$, $b\in\mathbb{F}_{2^m}^*$, and $m,n$ are integers satisfying $m|n$.

\begin{theorem}\label{tm1}
Let $f:\mathbb{F}_{2^n}\rightarrow \mathbb{F}_{2^m}$ be an arbitrary mapping, and $L:\mathbb{F}_{2^n}\rightarrow \mathbb{F}_{2^n}$ an $\mathbb{F}_{2^m}$-permutation of $\mathbb{F}_{2^n}$. Let $g:\mathbb{F}_{2^m}\rightarrow \mathbb{F}_{2^m}$ be a permutation of $\mathbb{F}_{2^m}$, and $A^{-1}(x)$ the compositional inverse of $A(x)$. If $\gamma\in \mathbb{F}_{2^n}^*$ is a $(b,A)$-linear translator of $f$, then
\begin{equation}\label{eq1}
\phi(x)=L(x)+L(\gamma)\left(g(f(x))+A^{-1}(\frac{f(x)}{b})\right)
\end{equation}
is a permutation polynomial over $\mathbb{F}_{2^n}$, and it compositional inverse is
$$\phi^{-1}(x)=L^{-1}(x)+\gamma A^{-1}\left(\frac{f(L^{-1}(x))}{b}+\frac{g^{-1}(A^{-1}(\frac{f(L^{-1}(x))}{b}))}{b}\right).$$
\end{theorem}
\textbf{Proof.}\quad We first consider the function $\varphi(x):\mathbb{F}_{2^n}\rightarrow \mathbb{F}_{2^n}$ as below:
$$\varphi(x)=x+\gamma\left(g(f(x))+A^{-1}(\frac{f(x)}{b})\right).$$

Putting $y=x+\gamma\left(g(f(x))+A^{-1}(\frac{f(x)}{b})\right)$, we derive that
\begin{align*}
f(y)&=f\left(x+\gamma\left(g(f(x))+A^{-1}(\frac{f(x)}{b})\right)\right)=f(x)+bA\left(g(f(x))+A^{-1}(\frac{f(x)}{b})\right)\\
     &=f(x)+bA(g(f(x)))+bA(A^{-1}(\frac{f(x)}{b}))=bA(g(f(x))),
\end{align*}
since $\gamma$ is a $(b,A)$-linear translator of $f$ and $A(x)$ is an additive polynomial. Further, we easily see that
$f(x)=g^{-1}(A^{-1}(\frac{f(y)}{b}))$, and
\begin{align*}
x&=y+\gamma\left(g(f(x))+A^{-1}(\frac{f(x)}{b})\right)
=y+\gamma A^{-1}\left(\frac{f(y)}{b}+\frac{g^{-1}(A^{-1}(\frac{f(y)}{b}))}{b}\right),
\end{align*}
implying that $\varphi(x)$ permutes $\mathbb{F}_{2^n}$ and its compositional inverse is
$$\varphi^{-1}(x)=x+\gamma A^{-1}\left(\frac{f(x)}{b}+\frac{g^{-1}(A^{-1}(\frac{f(x)}{b}))}{b}\right).$$

Now we point out that $\phi(x)=\varphi\circ L(x)$. Recall that $L$ is an $\mathbb{F}_{2^m}$-permutation of $\mathbb{F}_{2^n}$ such that $\phi(x)$ is a permutation of $\mathbb{F}_{2^n}$ as well. Moreover, $\varphi^{-1}(x)$ gives the sought for result:
$$\phi^{-1}(x)=L^{-1}(x)+\gamma A^{-1}\left(\frac{f(L^{-1}(x))}{b}+\frac{g^{-1}(A^{-1}(\frac{f(L^{-1}(x))}{b}))}{b}\right).$$ $\hfill\square$

\begin{corollary}
Let $f:\mathbb{F}_{2^n}\rightarrow\mathbb{F}_{2^m}$ be an arbitrary mapping, and $L:\mathbb{F}_{2^n}\rightarrow \mathbb{F}_{2^n}$ an $\mathbb{F}_{2^m}$-permutation of $\mathbb{F}_{2^n}$. Let $g:\mathbb{F}_{2^m}\rightarrow \mathbb{F}_{2^m}$ be a permutation. Suppose that $\gamma_1,\gamma_2,\gamma_3\in \mathbb{F}_{2^n}^*$ satisfying $\gamma_1+\gamma_2+\gamma_3\neq0$ are all pairwise distinct $(b,A)$-linear translators of $f$. Let $\rho(x)=g(f(x))+A^{-1}(\frac{f(x)}{b})$, $\tilde{\rho}=A^{-1}\left(\frac{f(x)}{b}+\frac{g^{-1}(A^{-1}(\frac{f(x)}{b}))}{b}\right)$, where $A^{-1}(x)$ is the compositional inverse of $A(x)$. Then

\begin{align*}
H(x,y)=&Tr(xL(y))+Tr(L(\gamma_1)x\rho(y))Tr(L(\gamma_2)x\rho(y))\\
   &+Tr(L(\gamma_1)x\rho(y))Tr(L(\gamma_3)x\rho(y))+Tr(L(\gamma_2)x\rho(y))Tr(L(\gamma_3)x\rho(y))
\end{align*}
is bent. Moreover, its dual function $\tilde{H}$ is
\begin{align*}
\tilde{H}(x,y)=&Tr(yL^{-1}(x))+Tr(\gamma_1y\tilde{\rho}(L^{-1}(x)))Tr(\gamma_2y\tilde{\rho}(L^{-1}(x)))\\
   &+Tr(\gamma_1y\tilde{\rho}(L^{-1}(x)))Tr(\gamma_3y\tilde{\rho}(L^{-1}(x)))
   +Tr(\gamma_2y\tilde{\rho}(L^{-1}(x)))Tr(\gamma_3y\tilde{\rho}(L^{-1}(x))).
\end{align*}

\end{corollary}
\textbf{Proof.}\quad  Define $\phi_i(x)=L(x)+L(\gamma_i)\left(g(f(x))+A^{-1}(\frac{f(x)}{b})\right)$ as expression (\ref{eq1}), where $i\in\{1,2,3\}$, indeed, $\phi_i(x)=L(x)+L(\gamma_i)\rho(x)$. The rest follows immediately by a highly similar argument to  Theorem 1 in \cite{Mesnager2017} or Theorem 6 in \cite{Cepak2019}.   $\hfill\square$
\begin{remark}
Mesnager \cite{Mesnager2014} provided a method for constructing bent functions and computed their duals. One showed that if $\phi_1,\phi_2,\phi_3$ are three pairwise distinct permutations of $\mathbb{F}_{2^n}$ such that $\psi=\phi_1+\phi_2+\phi_3$ is a permutation of $\mathbb{F}_{2^n}$ and $\psi^{-1}=\phi^{-1}_1+\phi^{-1}_2+\phi^{-1}_3$. Then $H(x,y)=Tr(x\phi_1(y))Tr(x\phi_2(y))+Tr(x\phi_1(y))Tr(x\phi_3(y))
+Tr(x\phi_2(y))Tr(x\phi_3(y))$ is bent over $\mathbb{F}_{2^n}\times \mathbb{F}_{2^n}$. Furthermore, the dual of $H$ is given by $\tilde{H}(x,y)=Tr(y\phi^{-1}_1(x))Tr(y\phi^{-1}_2(x))$ $+Tr(y\phi^{-1}_1(x))Tr(y\phi^{-1}_3(x))
+Tr(y\phi^{-1}_2(x))Tr(y\phi^{-1}_3(x))$.
\end{remark}

\medskip


 Using the classic definition of linear translators, Qin and  Yan  \cite[Theorem 2.1]{qin2020} in 2020 constructed a class of permutation polynomials of the shape $F(x)=x+\sum\limits_{i=1}\limits^{m}{\gamma_ih_i(f_i(x))}$. Here we show an analogous result with minor modifications. In addition, we also give their compositional inverses.

\begin{theorem}\label{tm3.1}
Let $f_1, f_2,\cdots, f_m:\mathbb{F}_{q^n}\rightarrow \mathbb{F}_q$, $h_1, h_2, \cdots, h_m:\mathbb{F}_q\rightarrow \mathbb{F}_q$. 
 Let $b_1, b_2, \cdots, b_m\in \mathbb{F}_q$, and $\gamma_1, \gamma_2, \cdots, \gamma_m\in \mathbb{F}_{q^n}^*$ be linearly independent over $\mathbb{F}_q$. If $\gamma_i$ is a $(b_i,A_i)$-linear translator of $f_i$ for $i\in\{1, 2, \cdots, m\}$, and a $(0,A_j)$-linear translator of $f_j$ when $j\neq i$, then
$$F(x)=x+\sum\limits_{i=1}\limits^{m}{\gamma_ih_i(f_i(x))}$$
is a permutation polynomial over $\mathbb{F}_{q^n}$ if and only if $g_i(x)=x+b_iA_i(h_i(x))$ is a permutation polynomial over $\mathbb{F}_q$, where $i\in \{1,2,\cdots, m\}$. Moreover, the compositional inverse of $F(x)$ is given by
 $$F^{-1}(x)=x-\sum\limits_{i=1}\limits^{m}{\gamma_ih_i(g_i^{-1}(f_i(x)))}.$$
\end{theorem}
\textbf{Proof.}\quad Take any $\alpha, \beta\in \mathbb{F}_{q^n}$ such that $F(\alpha)=F(\beta)$. Then
\begin{equation}\label{eq3.1}
\alpha+\sum\limits_{i=1}\limits^{m}{\gamma_ih_i(f_i(\alpha))}=\beta+\sum\limits_{i=1}\limits^{m}{\gamma_ih_i(f_i(\beta))},
\end{equation}
and hence
$$\alpha=\beta+\sum\limits_{i=1}\limits^{m}{\gamma_i\left(h_i(f_i(\beta))-h_i(f_i(\alpha))\right)}.$$
For the sake of convenience, we denote $h_i(f_i(\beta))-h_i(f_i(\alpha))=c_i\in \mathbb{F}_q$,
and thus the above equation can be rewritten as $\alpha=\beta+\sum\limits_{j=1}\limits^{m}{\gamma_jc_j}$. Therefore, Equation (\ref{eq3.1}) is equivalent to
$$\beta+\sum\limits_{j=1}\limits^{m}{\gamma_jc_j}+\sum\limits_{i=1}\limits^{m}{\gamma_ih_i(f_i(\beta+\sum\limits_{j=1}\limits^{m}{\gamma_jc_j}))}=\beta+\sum\limits_{i=1}\limits^{m}{\gamma_ih_i(f_i(\beta))}.$$
Tidying up both sides of the equation yields

\begin{equation*}
\sum\limits_{i=1}\limits^{m}\left(h_i(f_i(\beta+\sum\limits_{j=1}\limits^{m}\gamma_jc_j))-h_i(f_i(\beta))+c_i\right)\gamma_i=0,
\end{equation*}
which implies that for $1\leq i\leq m$,
$$h_i(f_i( \beta+\sum\limits_{j=1}\limits^{m}{\gamma_jc_j}))-h_i(f_i(\beta))+c_i=0,$$
since $\gamma_1,\gamma_2,\cdots, \gamma_m$ are linearly independent over $\mathbb{F}_q$.
Because $\gamma_i$ is a $(b_i,A_i)$-linear translator of $f_i$, and a $(0,A_j)$-linear translator of $f_j$ for $j\neq i$,
applying Proposition \ref{pro2.2}, we reformulate the equation above as
\begin{equation}\label{eq3.2}
h_i(f_i(\beta)+b_iA_i(c_i))-h_i(f_i(\beta))+c_i=0,
\end{equation}
where $i$ is an integer with $1\leq i\leq m$. Note that $F(x)$ is a permutation polynomial of $\mathbb{F}_{q^n}$
if and only if for all $i$, $c_i=0$. Clearly, the conclusion has been proven if $b_1=b_2=\cdots=b_m=0$.
Thus, we next consider the case in which some $b_i\neq 0$ for $1\leq i \leq m$.  Without loss of generality,
assume that the first $l$ $b_i$ are not zero and that the others are all zero. Then it follows from Equation (\ref{eq3.2}) that $c_{l+1}=\cdots=c_m=0$, and
$$A_i(h_i(f_i(\beta)+b_iA_i(c_i))+c_i)=A_i(h_i(f_i(\beta))),~~~ {\rm for} ~1\leq i\leq l,$$
which is equivalent to
$$A_i(h_i(f_i(\beta)+b_iA_i(c_i)))+b_i^{-1}(b_iA_i(c_i)+f_i(\beta))=b_i^{-1}f_i(\beta)+A_i(h_i(f_i(\beta))),$$
since $A_i(x)$ is an additive permutation over $\mathbb{F}_q$ and $b_i\neq 0$ for $1\leq i\leq l$. By a simple transformation, we also have
$$b_iA_i(h_i(f_i(\beta)+b_iA_i(c_i)))+f_i(\beta)+b_iA_i(c_i)=f_i(\beta)+b_iA_i(h_i(f_i(\beta))),~~{\rm for} ~1\leq i\leq l,$$
which implies that $g_i(f_i(\beta)+b_iA_i(c_i))=g_i(f_i(\beta))$, where $1\leq i\leq l$. We claim now that
for $1\leq i\leq l$, $c_i=0$ if and only if $g_i(x)$ is a permutation of $\mathbb{F}_q$ since $b_i\neq 0$.
By what we have already shown, we assert that $F(x)$ is a permutation of $\mathbb{F}_{q^n}$ if and only if
$g_i(x)$ is a permutation of $\mathbb{F}_q$.

Now it remains to compute the compositional inverse of $F(x)$. Obviously, the information provided above yields the following commutative diagram, and thus it follows from \cite[Theorem 3.1]{yuan2022-AGWPP} that
$F^{-1}(x)=x-\sum\limits_{i=1}\limits^{m}{\gamma_ih_i(g_i^{-1}(f_i(x)))}.$

\medskip

\centerline{\xymatrix@C=3em@R=2.5em{
\mathbb{F}_{q^n} \ar[r]^{F(x)} \ar[d]_{f_i} & \mathbb{F}_{q^n} \ar[r]^{F^{-1}(x)} \ar[d]^{f_i} & \mathbb{F}_{q^n} \ar[d]^{f_i} \\
\mathbb{F}_q \ar[r]_{g_i(x)} & \mathbb{F}_q\ar[r]_{g_i^{-1}(x)} & \mathbb{F}_q
}}

\medskip

For more clarity, we have $F(F^{-1}(x))=x$ since $F^{-1}(x)$ is the compositional inverse of $F(x)$, and so
$F(x)\circ F^{-1}(x)=(x+\sum\limits_{i=1}\limits^{m}{\gamma_ih_i(f_i(x))})\circ F^{-1}(x),$ that is
\begin{align*}
x&=F^{-1}(x)+\sum\limits_{i=1}\limits^{m}{\gamma_ih_i(f_i(F^{-1}(x)))}\\
   &=F^{-1}(x)+\sum\limits_{i=1}\limits^{m}{\gamma_ih_i(g_i^{-1}(f_i(x)))}.
\end{align*}
Consequently, $F^{-1}(x)=x-\sum\limits_{i=1}\limits^{m}{\gamma_ih_i(g_i^{-1}(f_i(x)))}.$  We are done. $\hfill\square$


\medskip

In \cite{Charpin2016}, Charpin, Mesnager and Sarkar utilized $0$-linear translators to characterize some involutions over $\mathbb{F}_{2^n}$ shaped like $G(x)+\gamma f(x)$. Observing the expressions for $F$ and $F^{-1}$ in Theorem \ref{tm3.1}, a similar consequence can be obtained as follows.  

\begin{corollary}\label{cor3.2}
The hypotheses are those of Theorem \ref{tm3.1} with $b_i=0$ for $1\leq i\leq m$. If $char(\mathbb{F}_q)=2$, then $F(x)=x+\sum\limits_{i=1}\limits^{m}{\gamma_ih_i(f_i(x))}$ is an involution of $\mathbb{F}_{q^n}$.
\end{corollary}

\begin{corollary}\label{cor3.3}
Let $f_1, f_2,\cdots, f_m:\mathbb{F}_{q^n}\rightarrow \mathbb{F}_q$, $h_1, h_2, \cdots, h_m:\mathbb{F}_q\rightarrow \mathbb{F}_q$, and $L$ be an $\mathbb{F}_q$-linear permutation of $\mathbb{F}_{q^n}$.  
Suppose that $b_1, b_2, \cdots, b_m\in \mathbb{F}_q$, and that $\gamma_1, \gamma_2, \cdots, \gamma_m\in \mathbb{F}_{q^n}^*$ are linearly independent over $\mathbb{F}_q$. If $\gamma_i$ is a $(b_i,A_i)$-linear translator of $f_i$ for $i\in\{1, 2, \cdots, m\}$, and a $(0,A_j)$-linear translator of $f_j$ when $j\neq i$, then
$$G(x)=L(x)+\sum\limits_{i=1}\limits^{m}{L(\gamma_i)h_i(f_i(x))}$$
is a permutation polynomial over $\mathbb{F}_{q^n}$ if and only if $g_i(x)=x+b_iA_i(h_i(x))$ is a permutation polynomial over $\mathbb{F}_q$, where $i\in \{1,2,\cdots, m\}$. Moreover,
 $$G^{-1}(x)=L^{-1}\left(x-\sum\limits_{i=1}\limits^{m}{L(\gamma_i)h_i(g_i^{-1}(f_i(x)))}\right).$$
\end{corollary}
\textbf{Proof.}\quad
It is not hard to find that 
$G(x)=L(x)\circ F(x)$, where $F(x)$ is shown in Theorem \ref{tm3.1}. Therefore, it follows from Theorem \ref{tm3.1} and the fact that $L$ is an $\mathbb{F}_q$-linear permutation of  $\mathbb{F}_{q^n}$ that $G(x)$ is a permutation over $\mathbb{F}_{q^n}$ if and only if $g_i(x)=x+b_iA_i(h_i(x))$ is a permutation over $\mathbb{F}_q$, where $1\leq i\leq m$.

Clearly, by the following commutative diagram, we have $f_i(x)\circ G^{-1}(x)=g_i^{-1}(x)\circ f_i(x)$. Hence, we immediately obtain the compositional inverse of $G(x)$.

\medskip

\centerline{\xymatrix@C=3em@R=2.5em{
\mathbb{F}_{q^n} \ar[r]^{G(x)} \ar[d]_{f_i} & \mathbb{F}_{q^n} \ar[r]^{G^{-1}(x)} \ar[d]^{f_i} & \mathbb{F}_{q^n} \ar[d]^{f_i} \\
\mathbb{F}_q \ar[r]_{g_i(x)} & \mathbb{F}_q\ar[r]_{g_i^{-1}(x)} & \mathbb{F}_q
}}     $\hfill\square$

\begin{corollary}
The hypotheses are those of Corollary \ref{cor3.3} with $b_i=0$,where $1\leq i \leq m$. And let $F(x)$ be as in
Theorem \ref{tm3.1}. If $char(\mathbb{F}_q)=2$ and $L$ is an involution which satisfies $L \circ F=F\circ L$,
then $G$ is an involution and $L$ fixes $\gamma_i$ for $1\leq i\leq m$.
\end{corollary}
\textbf{Proof.}\quad First notice that $F(x)$ is an involution on $\mathbb{F}_{q^n}$, i.e., $F\circ F(x)=x$. The condition $b_i=0$ also implies that $g_i(x)=x$ is an identity, where $1\leq i\leq m$, so by Corollary \ref{cor3.3}, $G(x)$ is a permutation of $\mathbb{F}_{q^n}$. Consider the composite of $G(x)$ with itself, we have
$G\circ G(x)=G\circ L\circ F(x)=L\circ F\circ L\circ F(x)=L(F\circ L(F(x)))=L(L\circ F(F(x)))=L(L(x))=x$, which implies that $G^{-1}(x)=G(x)$, meaning that $G$ is an involution.

Finally, by comparing the expressions for $G$ and $G^{-1}$ in Corollary \ref{cor3.3}, one may show that $L(\gamma_i)=\gamma_i$ for any $1\leq i\leq m$.     $\hfill\square$

\begin{theorem}\label{th2}
 Assume that $f_1, f_2,\cdots,f_m$ are mapping from $\mathbb{F}_{q^n}$ to $\mathbb{F}_q$, and $h_1, h_2, \cdots,h_m$
 are permutation polynomials of $\mathbb{F}_q$. Let $L:\mathbb{F}_{q^n}\rightarrow \mathbb{F}_{q^n}$ be an $\mathbb{F}_q$-linear
 polynomial with $Ker L\cap ImL=\{0\}$, 
and let $\gamma_1,\gamma_2,\cdots,\gamma_m$ be a basis for $Ker L$ over $\mathbb{F}_q$. If $\gamma_i$ is a $(b_i, A_i)$-linear translator of $f_i$, in which $b_i\neq 0$, and a $(0, A_j)$-linear translator of $f_j$ when $j\neq i$, then
   $$F(x)=L(x)+\sum\limits_{i=1}\limits^{m}{\gamma_ih_i(f_i(x))}$$
is a permutation polynomial over $\mathbb{F}_{q^n}$. 
\end{theorem}
\textbf{Proof.}\quad For convenience, we use $ \Lambda_L$ to denote the set $Ker L\cap ImL$. Let $\alpha, \beta\in \mathbb{F}_{q^n}$ satisfy $F(\alpha)=F(\beta)$. We get that
 $$L(\alpha)+\sum\limits_{i=1}\limits^{m}{\gamma_ih_i(f_i(\alpha))}=L(\beta)+\sum\limits_{i=1}\limits^{m}{\gamma_ih_i(f_i(\beta))},$$
which can be equivalently expressed as
\begin{equation}\label{eq3.3}
L(\alpha-\beta)=\sum\limits_{i=1}\limits^{m}{(h_i(f_i(\beta))-h_i(f_i(\alpha)))\gamma_i}.
 \end{equation}
 We now claim that the value of the above equation lies in $\Lambda_L$. First we deduce that $\sum\limits_{i=1}\limits^{m}{(h_i(f_i(\beta))-h_i(f_i(\alpha)))\gamma_i}\in Ker L$ since $\gamma_1,\gamma_2,\cdots,\gamma_m$ is a basis of $Ker L$ over $\mathbb{F}_q$ and $h_i(f_i(\beta))-h_i(f_i(\alpha))\in \mathbb{F}_q$. On the other hand, $\sum\limits_{i=1}\limits^{m}{(h_i(f_i(\beta))-h_i(f_i(\alpha)))\gamma_i}\in Im L$ is a straightforward consequence from Equation (\ref{eq3.3}). Hence we arrive at our assertion.

Then by $\Lambda_L=\{0\}$,
 we have $\alpha-\beta\in Ker L$ is a linear combination $\alpha-\beta=\sum\limits_{j=1}\limits^{m}{a_j\gamma_j}$
 with  $a_j\in \mathbb{F}_q$ for $1\leq j \leq m$. Meanwhile, using the fact that $\gamma_1,\cdots,\gamma_m$ are
 linearly independent over $\mathbb{F}_q$, we obtain that $h_i(f_i(\beta))-h_i(f_i(\alpha))=0$, which yields that
  $f_i(\beta)-f_i(\alpha)=0$ since $h_i$ is a permutation of $\mathbb{F}_q$, where $1\leq i \leq m$.
 Combining these two aspects, we see for $1\leq i\leq m$, $$f_i(\beta)-f_i(\beta+\sum\limits_{j=1}\limits^{m}{a_j\gamma_j})=0,$$
and then based on the assumptions between $\gamma_i$ and $f_j$, for $1\leq i,j \leq m$,
we infer $b_iA(a_i)=0$ for $1\leq i\leq m$. Notice that $b_i\neq 0$ and $A_i(x)$ is an additive permutation, we have $a_i=0$ and, in consequence, $\alpha=\beta$, which indicates that $F(x)$ is a permutation of $\mathbb{F}_{q^n}$.    $\hfill\square$

\begin{corollary}
Let $p,m,b$ be integers with $p$ odd prime and $b\neq 0$. Let $f:\mathbb{F}_{p^{2m}}\rightarrow \mathbb{F}_{p^m}$ be an arbitrary mapping, $h$ a permutation of $\mathbb{F}_{p^m}$, and $\gamma\in \mathbb{F}_{p^{2m}}^*$ satisfies $\gamma+\gamma^{p^{m}}=0$. If $\gamma$ is a $(b,A)$-linear translator of $f$, then
 $$F(x)=x+x^{p^m}+\gamma h(f(x)) $$
is a permutation polynomial over $\mathbb{F}_{p^{2m}}$.
\end{corollary}
\textbf{Proof.}\quad To begin with, observe that for $Tr^{2m}_{m}(\gamma)=\gamma+\gamma^{p^{m}}=0$ we have $\gamma\in Ker(Tr^{2m}_{m}(x))$. Then $\gamma$ is a basis of $Ker (Tr^{2m}_{m}(x))$ over $\mathbb{F}_{p^{m}}$, since the dimension theorem of linear mapping indicates that $Ker (Tr^{2m}_{m}(x))$  as a $\mathbb{F}_{p^{m}}$-vector space is one dimension.

 Now we assert that $Ker (Tr^{2}_{m}(x))\cap Im (Tr^{2m}_{m}(x))=\{0\}$. Suppose, on the contrary,  that there exists $\xi \in  Ker (Tr^{2m}_{m}(x)) \cap  Im (Tr^{2m}_{m}(x))$
 with $\xi \neq 0$. Notice that $Tr^{2m}_{m}(x)$ maps $\mathbb{F}_{p^{2m}}$ onto $\mathbb{F}_{p^{m}}$,
 one shows that $Im (Tr^{2m}_{m}(x))$ $=\mathbb{F}_{p^{m}}$, and so $\xi\in {\mathbb{F}_{p^{m}}}$.
 Then $Tr^{2m}_{m}(\xi)=\xi+\xi^{p^{m}}=2\xi=0$, which contradicts the fact that $p$ is an odd prime. And the remainder immediately follows from Theorem \ref{th2}.  $\hfill\square$

In what follows, we consider the case that $\gamma_i$ is a $(b_{ij},x^{p^t})$-linear translator of $f_j$ for $1\leq i,j \leq m$.

\begin{theorem}
Let $q=p^k$. Assume that $f_1, f_2,\cdots,f_m:\mathbb{F}_{q^n} \rightarrow \mathbb{F}_q$, and $h_1, h_2, \cdots,h_m$ are permutations of $\mathbb{F}_q$. Let $L:\mathbb{F}_{q^n}\rightarrow \mathbb{F}_{q^n}$ be an $\mathbb{F}_q$-linear polynomial with $Ker L\cap ImL=\{0\}$, and $A(x)=x^{p^t}: \mathbb{F}_q\rightarrow \mathbb{F}_q$, where $t\leq k$. Let $\gamma_1,\gamma_2,\cdots,\gamma_m$ be a basis of $Ker L$ over $\mathbb{F}_q$. If $\gamma_i$ is a $(b_{ij}, A)$-linear translator of $f_j$ for $1\leq i,j \leq m$, then
$$F(x)=L(x)+\sum\limits_{i=1}\limits^{m}{\gamma_ih_i^{p^t}(f_i(x))}$$
permutes $\mathbb{F}_{q^n}$ if and only if rank$(B)=m$, where $B=(b_{ij}^{p^{k-t}})$.
\end{theorem}
\textbf{Proof.}\quad Let rank$(B)=m$ and $\alpha, \beta \in \mathbb{F}_{q^n}$ satisfy $F(\alpha)=F(\beta)$.
With a similar discussion prior as in the proof of Theorem \ref{th2}, we see that
$$f_i(\beta)-f_i(\beta+\sum\limits_{j=1}\limits^{m}{a_j\gamma_j})=0,~~{\rm for} ~1\leq i\leq m,$$
where $a_j\in \mathbb{F}_q$ are coefficients of the linear combination $\alpha=\beta+\sum\limits_{j=1}\limits^{m}{a_j\gamma_j}$. Since $\gamma_i$ is a $(b_{ij}, A)$-linear translator of $f_j, 1\leq i\leq m, 1\leq j \leq m$, we have
$$\sum\limits_{j=1}\limits^{m}{b_{ij}A(a_j)}=0,~~{\rm for} ~1\leq i\leq m.$$
Substituting  $A(x)=x^{p^t}$ into the above equation and then rasing both sides to the $p^{k-t}$-th power, we obtain
$$\sum\limits_{j=1}\limits^{m}{b_{ij}^{p^{k-t}}a_j}=0,~~{\rm for} ~1\leq i\leq m,$$
and we reformulate this in terms of matrix product, then it becomes
\begin{equation*}
\left[
\begin{array}{cccc}
 b_{11}^{p^{k-t}} & b_{12}^{p^{k-t}}& \cdots & b_{1m}^{p^{k-t}}\\
 b_{21}^{p^{k-t}} & b_{22}^{p^{k-t}}& \cdots & b_{2m}^{p^{k-t}}\\
  \vdots   &  \vdots & \ddots     &\vdots \\
 b_{m1}^{p^{k-t}} & b_{m2}^{p^{k-t}}& \cdots & b_{mm}^{p^{k-t}}
\end{array}
\right ]
\left[
\begin{array}{cccc}
 a_{1}\\
 a_{2}\\
 \vdots \\
 a_{m}
\end{array}
\right ]
=
\left[
\begin{array}{cccc}
 0\\
 0\\
 \vdots  \\
 0
\end{array}
\right ],
\end{equation*}
which implies that $(a_1,a_2,\cdots,a_m)$ is a solution of $BX=\boldsymbol{0}$,
where $B=(b_{ij}^{p^{k-t}}), X=(x_1,x_2,\cdots,x_m)^T$. Recall that rank$(B)=m$,
we conclude that $(a_1,a_2,\cdots,a_m)=(0,0,\cdots,0)$ is the unique solution of the system of equations above.
Hence we get $\alpha=\beta$, and so $F(x)$ is a permutation of $\mathbb{F}_{q^n}$.

Conversely, suppose $X=(a_1,a_2,\cdots,a_m)^T$ is a solution of $BX=\boldsymbol0$. Then we have
$\sum\limits_{j=1}\limits^{m}{b_{ij}a_j^{p^t}}=0$, i.e., $\sum\limits_{j=1}\limits^{m}{b_{ij}A(a_j)}=0$,
where $1\leq i\leq m.$ And since $\gamma_i$ is a $(b_{ij},A)$-linear translator of $f_j$ for $1\leq i,j\leq m$,
this may be viewed as
\begin{center} {$f_i(x+\sum\limits_{j=1}\limits^{m}{a_j\gamma_j})-f_i(x)=0$, ~~for $1\leq i\leq m$,}
\end{center}
where $x\in \mathbb{F}_{q^n}$ is arbitrary. Note that $h_1,h_2,\cdots,h_m$ are permutations of $\mathbb{F}_q$, we have
$$h_i(f_i(x+\sum\limits_{j=1}\limits^{m}{a_j\gamma_j}))=h_i(f_i(x)),~~{\rm for} ~1\leq i\leq m,$$
and therefore
$$h_i^{p^t}(f_i(x+\sum\limits_{j=1}\limits^{m}{a_j\gamma_j}))=h_i^{p^t}(f_i(x)),~~{\rm for} ~1\leq i\leq m.$$
Since $\gamma_1,\gamma_2,\cdots,\gamma_m$ are linearly independent over $\mathbb{F}_q$, we easily obtain
\begin{equation}\label{eq3.4}
\sum\limits_{i=1}\limits^{m}\left(h_i^{p^t}(f_i(x+\sum\limits_{j=1}\limits^{m}a_j\gamma_j))-h_i^{p^t}(f_i(x))\right)\gamma_i=0.
\end{equation}
Now $$L(x+\sum\limits_{j=1}\limits^{m}{a_j\gamma_j})=L(x)+\sum\limits_{j=1}\limits^{m}{a_jL(\gamma_j)=L(x)},$$
since $L$ is an $\mathbb{F}_q$-linear polynomial, along with $\gamma_j \in Ker L$ for $1\leq j\leq m$,
and so Equation (\ref{eq3.4}) is equivalent to
$$L(x+\sum\limits_{j=1}\limits^{m}{a_j\gamma_j})
+\sum\limits_{i=1}\limits^{m}{\gamma_ih_i^{p^t}(f_i(x+\sum\limits_{j=1}\limits^{m}{a_j\gamma_j}))}
=L(x)+\sum\limits_{i=1}\limits^{m}{\gamma_ih_i^{p^t}(f_i(x))},$$
that is,
$$F(x+\sum\limits_{j=1}\limits^{m}{a_j\gamma_j})=F(x).$$
Then we immediately deduce that $\sum\limits_{j=1}\limits^{m}{a_j\gamma_j}=0$ as $F(x)$ permutes $\mathbb{F}_{q^n}$.
 Further, it follows again from the $\mathbb{F}_q$-linear independence of $\{\gamma_1,\gamma_2,\cdots,\gamma_m\}$ that
 $a_1=a_2=\cdots=a_m=0$. This indicates that the system of equations $BX=\boldsymbol0$ has only zero solutions,
 and thus rank$(B)=m$.  $\hfill\square$
\begin{remark}
The idea of this conclusion comes from the work of Qin and Yan \cite[Theorem 2.4]{qin2020}, utilizing classical linear translators, they constructed a class of permutation polynomials shape like $F(x)=x+\sum\limits_{i=1}\limits^{m}\lambda_i\gamma_i f_i(x)$.
\end{remark}

%

\section*{Acknowledgements}
Pingzhi Yuan was supported by the National Natural Science Foundation of China (Grant No. 12171163) and Guangdong Basic and Applied Basic Research Foundation  (Grant No. 2024A1515010589).
Hongjian Li was supported by the Project of Guangdong University of Foreign Studies (Grant No. 2024RC063).


%
\section*{Declarations}
\begin{conflict of interest} {\rm There is no conflict of interest.}
\end{conflict of interest}

\end{document}